\magnification =1064
\hsize 16truecm
\vsize 23truecm
\hfuzz =5pt
\scrollmode
\def\qed{\hbox{\vrule height 7pt depth 0pt width 7pt}}
\def\cqfd{\hfill\penalty 500\kern 10pt\qed\medbreak}

\def \l{{\lambda}}
\def \a{{\alpha}}
\def \b{{\beta}}
\def \s{{\sigma}}

\def \o{{\omega}}
\def \O{{\Omega}}

\def \R{{\bf R}}
\def \Z{{\bf Z}}

\def \Q{{\bf Q}}\def \m{{\mu}}

\def \E{{\bf E}\,}

\def \A{{\cal A}}

 \def \P{{\bf P}}

\def \qq{{\qquad}}
 \def \p{{\varphi}}
\def \noi{{\noindent}}

\def \e{{\varepsilon}}
 \def \t{{\theta}}
\def \T{{\bf T}}

\def\YY{ {{\bf{Y}}}} 
\def\th{ {{ \underline{t}}}}

\def \lt {{\hbox{\vrule height 8pt depth 2pt width 0,9pt}\kern 1,1pt}}
\def \rt {{\kern 1,1pt\hbox{\vrule height 8pt depth 2pt width 0,9pt}}}
 
\font\ph=cmcsc10  at  10 pt
\font\iit =cmmi7 at 7pt

\font\ph=cmcsc10  at  10 pt

\font\gum= cmcsc10 at 11 pt
\font\gem= cmbx10 at 10 pt
\font\gum= cmbx10 at 11 pt

\font\ph=cmcsc10  at  10 pt

\font\iit =cmmi7 at 7pt
 \vskip 10mm
\font\gum= cmbx10 at 11 pt

 at 10,7 pt
\font\gem= cmbx10 at 10 pt

\vskip 2cm

\rm
\def\ddate {\ifcase\month\or January\or
February\or March\or April\or May\or June\or July\or
August\or September\or October\or November\or December\fi\ {\the\day},
{\the\year}}

\rm
 \centerline{\gum On   localization   in Kronecker's diophantine theorem}
      \vskip 0,5 cm
  \centerline{ Michel {\ph Weber}}
\vskip 0,5 cm
  {\leftskip =1,8cm \rightskip=1,8cm   \noindent {\bf Abstract:\rm\  Using a probabilistic approach, we extend for general $\Q$-linearly
independent sequences a   result of T\'uran  concerning   the sequence $ ( \log p_\ell)$, $p_\ell$ being the $\ell$-th prime. For
instance  let
$   
\l_1,
\l_2,
\ldots  $   be linearly independent over $\Q$.  We prove that there  
exists a   constant $C_0$ such that for any positive integers    $N $ and $\o $,   if 
 $T> \big(   {4\o\over C_0} \sqrt{\log {N\o\over C_0} } 
        \big)^{N} / \Xi   $,   where 
 \vskip -8pt$$\Xi= \min_{{ u_k \,{\rm integers} \atop |u_k|\le 6 \o 
\log {(N\o/ C_0)}}\atop  |u_1\l_1+\ldots+ u_N\l_N |\not=0 } \big |\sum_{1\le k\le N }    \l_k 
 u_k \big|$$  
\vskip -4pt \noi  then to     any reals $d,  
\b_1,\ldots, 
\b_N  $,     corresponds a real $t\in [d,d+T]$  such that}  
 $ \sup_{j=1}^N\, \lt  t\l_j-\b_j \rt\le {1/ \o}$.      
   \par } 
  
\footnote {}{ \sevenrm  \hskip -2pt [Loc]Kronecker \ {\iit Date}\sevenrm : on \ddate \ \par \vskip -2pt   {\iit AMS} \ {\iit Subject}\
{\iit Classification} 2000:  Primary 11K60,
 Secondary 60G50, 11J25. \par \vskip -2pt
  {\iit Keywords}:     Kronecker's theorem, localization, Tur\'an's theorem, diophantine approximation, small deviations, weighted i.i.d.
sums. 
  \par}

 \vskip 0,7cm\noi

\medskip\noi {\gum 1. Introduction and main result}
 \medskip\noi  The well-known theorem of Kronecker on Diophantine approximation asserts that if $   \l_1, \l_2, \ldots ,\l_N $   are
linearly independent over $\Q$, then for any given real numbers $     \a_1, \a_2, \ldots , \a_N $  and any $\e >0$, there exists a real
number $t$ such that $$ \sup_{j=1}^N\, \lt  t\l_j-\a_j \rt\le  \e ,  \eqno(1.1) $$ 
where 
$\lt x
\rt$ denotes the distance of
$x$ to
$\Z$, i.e.  $  {\lt x
\rt=
\min_{\nu\in \Z} |x-\nu|}$. 

\smallskip A quantitative form of Kronecker's theorem was given by Bacon [B], who proved that if   $   \l_1, \l_2,
\ldots
,\l_N $  are reals numbers satisfying for some $M\ge 1$ $$\line{$ \qq   \quad
\left\{  
 \matrix{ u_1\l_1+\ldots+u_N\l_N =0 \cr   
  |u_1|+\ldots+|u_N|\le M,     \    \hbox{\rm $u_k$ integers}  \cr}  \right.  \qq \Longrightarrow\qq    
   u_1=u_2=\ldots=u_N=0,       \hfill$}  \eqno(1.2) $$
then for any   real numbers $     \a_1, \a_2, \ldots , \a_N $, there exists a real
number $t$ such that $$ \sup_{j=1}^N\, \lt  t\l_j-\a_j \rt\le {c(N)\over M}, \qq \quad c(N)= {1\over 2}(N-1)^{3/2}\left({125\over
48}\right)^{(N^3-N )/12} .  
 \eqno(1.3)$$

Recently  Chen [C1] considerably improved   this result, showing that there exists a real
number $t$ such that 
 $$ \sum_{n=1}^N \lt t\l_n-\a_n\rt \le  {\pi^2\over 16}{N \over   (M+1)^2} . \eqno(1.4)$$ He also considered the case when  $   \l_1,
\l_2,
\ldots
\l_N
$ and
$    
\a_1,
\a_2, \ldots ,
\a_N
$ are real numbers such that for some $M\ge 1$ 
$$\line{$ \qq   \quad
\left\{  
 \matrix{ u_1\l_1+\ldots+u_N\l_N  \ \hbox{\rm is an integer} \cr   
  |u_1|+\ldots+|u_N|\le M,     \    \hbox{\rm $u_k$ integers}  \cr}  \right.  \qq \Longrightarrow\qq    
  u_1\a_1+\ldots+u_N\a_N \ \hbox{\rm is an integer}.       \hfill$}   \eqno(1.5)$$
 
 No indication is however  given on the range of $t$, and in [C1] it was claimed that no estimate for $t$ exists in general.  We refer
to   [Tu] (see also [Ti]) for more information about this important facet of Kronecker's theorem.
  The object of
this work is to provide a simple  estimate for
$t$.        \medskip\noi {\gem Theorem 1.} {\it   There   exists a   constant $C_0$ such that for any positive integers    $N $,
$\o $,     if  
$\l_1,
\l_2,
\ldots,
\l_N 
$ are reals satisfying  
 $$\line{$ \qq   \quad
\left\{  
 \matrix{ u_1\l_1+\ldots+u_N\l_N =0 \cr   
 \displaystyle{\max_{1\le \ell\le N}}
|u_\ell|\le 6 \o  \log {N\o\over C_0} ,     \    \hbox{\rm $u_k$ integers}  \cr}  \right.  \qq \Longrightarrow\qq    
   u_1=u_2=\ldots=u_N=0,       \hfill$}\eqno(1.6)   $$
 
  if    
$$T    >  {3 \over \pi \Xi} \bigg( { 2\sqrt 3\,\o  \over   
      C_0} \sqrt{ \log \ {N\o\over C_0}   }\bigg)^{N}  \qq {\it where}\qq\Xi=\Xi(N,\o):= \min_{{ u_k \,{\rm integers} \atop |u_k|\le
  6\o  \log {(N\o/ C_0)}}\atop  |u_1\l_1+\ldots+ u_N\l_N |\not=0 } \big |\sum_{1\le k\le N }    \l_k 
 u_k \big|, $$
  then to     any reals $d,  
\b_1,\ldots, 
\b_N  $     corresponds a real $t\in [d,d+T]$  such that}  
$$ \sup_{j=1}^N\, \lt  t\l_j-\b_j \rt\le {1\over \o}.  \eqno(1.7)  $$

  \bigskip  
    When $    \l_1, \l_2, \ldots, \l_N 
$ are linearly independent over $\Q$, condition (1.6) is trivially satisfied, and so the theorem applies.  
  In the  case $\l_\ell= \log p_\ell$, $p_\ell$ being the $\ell$-th prime,  $\ell =1,\ldots, N$, Tur\'an    ([Tu], Lemma p.313)
proved that the conclusion above is satisfied with $T= e^{17\o N\log^2 N}$ if $N$ is large enough, and $  4\le \o\le N$. It is
possible to estimate
$\Xi$ from below. More precisely,  given any real
$\e>0$, there exists and integer
$N(\e)$ depending on $\e$  only, such that for $N\ge N_\e$. 
$$ \Xi  \ge      e^{ - (1+\e) \o  N\log {(N\o/ C_0)}  \log N} $$
 From this and Theorem 1, we deduce  the similar estimate: if $N$ is large enough, and $\o$ is any positive integer, one can take  $T   
>   e^{   (1+2\e)
\o  N\log {(N\o/ C_0)} 
\log N  }   $.
\smallskip  The proof of
Theorem 1 is inspired from T\'uran's proof of the aforementionned particular case. But we also introduced an important probability
structure allowing us to tacle the general case.
\medskip
Let us make some further remarks.    By Theorem 1,  we can take $T=    \big(   {4\o\over C_0} \sqrt{\log {N\o\over C_0} } 
        \big)^{N} / \Xi  $. 
Then $   { \o\over C_0} \sqrt{\log {N \o\over C_0} } ={ (T \Xi)^{1/N}/ 4}   $. Let $w={ \o/ C_0}$, $\Theta={ (T \Xi)^{1/N}/ 4}$. As $   w
\sqrt{\log Nw } =\Theta  $, we get
$${\Theta\over \sqrt{\log N\Theta}} ={w\sqrt{\log Nw }\over \sqrt{\log N(w\sqrt{\log w })}}\le w ={ \o\over C_0}.$$
  Since by (1.7),  to     any reals $d,  
\b_1,\ldots, 
\b_N  $, corresponds a real $t\in [d,d+T]$  such that  $ \sup_{j=1}^N\, \lt  t\l_j-\b_j \rt\le {1/ \o}  $, we are free to choose $d=T/2$. 
Then 
$$  \sup_{j=1}^N\, \lt  t\l_j-\b_j \rt\le {1\over \o} \le  {   \sqrt{\log N\Theta} \over C_0\Theta}= {   \sqrt{\log { N { (T \Xi)^{1/N}/
4}}} \over C_0{ { (T
\Xi)^{1/N}/ 4}}}\le 4 {   \sqrt{\log     t  N { (  \Xi)^{1/N} } } \over     C_0     { ( t \Xi)^{1/N} } } , $$
or 
$$ {      { (t  \Xi)^{1/N} }\over  \sqrt{\log     t  N { (  \Xi)^{1/N} } }      }\sup_{j=1}^N\, \lt  t\l_j-\b_j \rt\le {4 \over C_0}.$$
And this holds  for infinitely many $t$. We deduce  
$$\liminf_{t\to \infty}{  t^{1/N}  \over  \sqrt{ \log      t } }  \cdot \displaystyle{\sup_{j=1}^N\, \lt  t\l_j-\b_j \rt  } 
<\infty.$$
  In particular for $\l_1$ irrational  
 $$\liminf_{t\to \infty}{  t\,   \lt  t\l_1  \rt  \over   \sqrt{ \log      t }  }    <\infty.$$
and   if $\l_1$,   $\l_2$ and $\l_1/\l_2$  
are irrationals  
$$\liminf_{t\to \infty}{  t\,   \lt  t\l_1  \rt\ \lt  t\l_2  \rt \over   \log     t    }     <\infty.$$

The well-known Littlewood's conjecture       (see [M] p. 202), however  states  that for any   $\l_1$,   $\l_2$  
irrationals 
 $\displaystyle{\liminf_{t\to \infty}   t\,   \lt  t\l_1  \rt\ \lt  t\l_2  \rt       =0}$. 

 Finally, applications of Theorem 1   to supremums of Dirichlet polynomials and more general polynomials are given at the end of Section
3.

  \bigskip

\noi{\gum  2. Some probabilistic preliminaries}

   \medskip\noi  
Let $e(x)=e^{2i\pi x}$.  Let $m$ be a positive  integer. Let $(\O,\A,\P)$ be a probability space,  and let  $X$ be a discrete   random
variable with law defined by:
 $$\line{$  \qq\qq\qq\qq\qq\qq  \P\{X=n\} =
 \left\{  
 \matrix{ {m-|n|\over m^2}         &  \qq {\rm if} \ 0\le |n|<m, \cr   
 0          &   {\rm if} \ |n|\ge m.  \cr}   \right.   \hfill$}$$
  Then
 $\E X= 0$, $\s^2:  =\E X^2 = {(m^2 -1)/ 6}$, and the characteristic function $\p_X(t)=\E\, e(tX) $      satisfies 
$$\p_X(t)=\sum_{0\le |n|<m} \P\{X=n\}e(tn)  ={1\over m^2} \sum_{0\le |n|<m} (m-|n|)  e(tn) =m^{-2}\cdot
|A_m(e(t))|^2,$$
 where    $ A_m(z) =1+z+\ldots+z^{m-1} $. Indeed  we have
$$|A_m(z)|^2 =\sum_{j=0}^{m-1}\sum_{\ell=0}^{m-1} z^{j-\ell}=\sum_{n=-m+1}^{m-1} z^n\#\{0\le j,\ell<m:  j-\ell=n\}=\sum_{0\le |n|<m}
\big(m-|n|)z^{  n }.  $$
\medskip {\gem Remark 1.}  
--- We have $\p_X(t)=  (2 \pi/m) F_m(2\pi t)   $, where $F_m$ is the Fej\'er kernel
$$  F_m(u) ={1\over 2m\pi} \Big({\sin mu/2\over
\sin u/2}\Big)^2 = {1\over m}\sum_{k=0}^{m-1} D_k(u),\qq\qq D_m(u)={1\over 2\pi}\sum_{|k|\le m}e^{-iku}, $$
$D_m$ being the Dirichlet kernel.
 Now let $X_1,\ldots, X_k$ be independent copies of $X$. Put $S_k=X_1+\ldots+ X_k$, and consider its characteristic function  $
\p_{S_k}(t)=\E\, e(tS_k)$. Basic properties of independent random variables imply  
$$ \p_{S_k}(t) =\p^k_X(t)    =\sum_{0\le |\nu|\le (m-1)k} \P\{S_k=\nu\} e(t\nu ) =m^{-2k}\cdot |A_m(e(t))|^{2k} . \eqno(2.1)$$
  By the local limit theorem [P] p.187
 $$ \sup_{\nu}\Big|\s \sqrt k\P\{S_k=\nu\}-{1\over \sqrt{2\pi}}e^{-{\nu^2\over 2\s^2k}}\Big|=
o(1)\qq k\to\infty.$$ Thereby
$$\P\{S_k=\nu\}= {1\over \sqrt{ \pi k(m^2-1)/3}}e^{-{3\nu^2\over  (m^2-1)k}}+o(1){1\over \sqrt{k(m^2-1)/6}},$$
and in particular for each $m$, as $k$ tends to infinity 
$$  \P\{S_k=0\}=m^{-2k}\int_0^1 |A_m(e(t))|^{2k}dt=\int_0^1 |{\sin  { \pi  mt  }
\over m
\sin{ \pi  t } }|^{2k}dt =
\big({3\over\pi}\big)^{1/2}{1\over
m\sqrt{
  k }}(1+o(1)) .\eqno(2.2)$$  

  When $m$  and $k$ vary simultaneously,   some useful estimates are also at disposal ([Tu]). 
  For $k$ large, and any positive integer $m$  
$$ \int^{1 }_0 \big({  \sin    {  \pi  mt  }
\over
 \sin  {   \pi t } }\big)^{2 k}dt  \ge
C{m^{2k-1}\over \sqrt k},  $$
where $C$ is an absolute constant.  
 Indeed, with the variable change $t=u/m\pi$
$$\eqalign{ \int_0^1 |{\sin  { \pi  mt  } \over
\sin{ \pi  t } }|^{2k}dt&={1\over m\pi} \int_0^{m\pi} |{\sin  {    u  } \over
\sin{    u/m } }|^{2k}du \ge {1\over m\pi}\int_0^{m\pi} |{\sin  {    u  } \over
      (u/m)  }|^{2k}du 
\cr & \ge {m^{2k-1}\over  \pi}  \int_0^{m\pi} |{\sin  {    u  } \over
       u   }|^{2k}du
 \ge {m^{2k-1}\over  \pi}\int_0^{k^{-1/2} } |{\sin  {    u  } \over
       u   }|^{2k}du\cr &\ge {m^{2k-1}\over  \pi}\int_0^{k^{-1/2} } |1-{u^2\over 5}|^{2k}du
  = {m^{2k-1}\over  \pi}\int_0^{k^{-1/2} } e^{{2k}\log (1-{u^2\over 5} )}  du
\cr &\ge {m^{2k-1}\over  \pi}\int_0^{k^{-1/2} } e^{- { ku^2\over  5} }  du\ge
C{m^{2k-1}\over \sqrt k}. \cr} $$

And so there exist constants $k_0$, $C_0>0$ such that for $k\ge k_0$ and any positive integer $m$
 $$  \P\{S_k=0\}\ge {C_0\over m\sqrt k} .\eqno(2.3)$$ 
Further, we may and do assume $C_0< 1/4$.  Conversely, notice that for any $m>m_0$, and any positive $k$ 
  $$ \int_0^1 |{\sin  { \pi  mt  } \over
\sin{ \pi  t } }|^{2k}dt \le 2^{2k+1} m^{4k^2 /(2k+1)}=C\cdot m^{2k-1+ 1 /(2k+1)}. $$

\bigskip

 Finally, as $\sin \pi   t\ge ({2/ \pi })\pi t=2t$, for $0\le t\le 1/2$, we have
$$     \p_{S_k}(t)    =\Big({\sin  \pi m
t\over
m\sin \pi   t}\Big)^{2k}  \le \Big({ 1\over
2m    \lt t \rt}\Big)^{2k}  \wedge 1.
\eqno(2.4) $$
\medskip
 \medskip\noi {\gum  3. Proof of Theorem 1.}
 \medskip \noi  Let $    
\b_1,\ldots, 
\b_N$ be given reals.  
Let $Y_1,\ldots, Y_N$ be independent copies of $S_k$. Consider the random vector $\YY=(Y_1,\ldots, Y_N) $ and let 
$\underline{  \b}=( 
\b_1,\ldots, 
\b_N) $, $\th= (t\l_1-
\b_1,\ldots,t\l_N-
\b_N)$. Put
$${\bf \Upsilon}(t, \underline{  \b}):  =\E\,  e(\langle \th, \YY\rangle)  =\E\,  e  
\Big(t\sum_{\ell=1}^N \l_\ell Y_\ell-\sum_{\ell=1}^N\b_\ell Y_\ell\Big) =\prod_{\ell=1}^N \p_{S_k}(t\l_\ell-
\b_\ell) , $$
and for $j=1,\ldots, N$ 
$${\bf \Upsilon}_j(t, \underline{  \b}):  =\E\,  e \Big(t\sum_{1\le \ell\le N \atop \ell \not= j } \l_\ell Y_\ell-\sum_{1\le
\ell\le N \atop \ell \not= j }\b_\ell Y_\ell\Big) =\prod_{1\le \ell\le N \atop \ell \not= j } \p_{S_k}(t\l_\ell-
\b_\ell) . $$
Let $\o\ge 1$.  Let $d    
 $ be  another given real  and let $T>0$. Suppose that  to any $t\in [d,d+T]$, corresponds an indice
$j=j_t\in \{1,\ldots,N\}$, such that  
$$  \lt  t\l_j-\b_j \rt>{1/\o}.  \eqno(3.1) $$  
We will show that this can happen only if $T$ is not too large. By (2.4)
 
$$     \p_{S_k}(t\l_j-
\b_j) \le \Big({ 1\over
2m    \lt t\l_j-
\b_j \rt}\Big)^{2k}   \le \big({\o\over
2m     }\big)^{2k}  .
 $$ 
and   so
$${\bf \Upsilon}(t, \underline{  \b})\le  \big({\o\over
2m     }\big)^{2k} \sum_{j=1}^N\chi\{ j_t=j\} {\bf \Upsilon}_j(t, \underline{  \b}).   $$
 Integrating this inequality over  $  [d,d+T]$ yields
$$\eqalign{\int_{d}^{d+T}{\bf \Upsilon}(t, \underline{\b})dt &\le   \big({\o\over
2m     }\big)^{2k}  \sum_{j=1}^N\int_{d}^{d+T} \chi\{ j_t=j\} {\bf \Upsilon}_j(t, \underline{  \b}) dt  \le   \big({\o\over
2m     }\big)^{2k} \sum_{j=1}^N\int_{d}^{d+T}    {\bf \Upsilon}_j(t, \underline{  \b}) dt. \cr}\eqno(3.2) $$
 
   But 
 $$\eqalign{ {\bf \Upsilon}_j(t, \underline{  \b})&= \prod_{1\le \ell\le
N \atop j\not=
\ell} \p_{S_k}(t\l_\ell-
\b_\ell)  = \prod_{1\le \ell\le
N \atop \ell\not=
j} \bigg(\sum_{0\le |\nu|\le (m-1)k} \P\{S_k=\nu\} e\big( (t\l_\ell-
\b_\ell)
\nu\big) \bigg)  \cr &=\P\{S_k=0\}^{N-1}\cr &+\sum_{0<  {\sup_{ \ell\not=j}}|\nu_\ell|\le (m-1)k  } \Big(\prod_{1\le \ell\le
N \atop j\not=
\ell} \P\{S_k=\nu_\ell\} \Big)e(- \sum_{1\le \ell\le
N \atop j\not=
\ell}  \b_\ell 
\nu_\ell)\cdot e( t\sum_{1\le \ell\le
N \atop j\not=
\ell}  \l_\ell 
\nu_\ell).\cr}\eqno(3.3)  $$
Let $C_0$ be  the constant from (2.3). As $C_0< 1/4$, it follows that
    $  N\o 
    >   4C_0 $. Choose   
 $$m= 2\o    ,\qq k= \inf\Big \{j\ge 1:  {N\o\over
 C_0     }    \le  {4^{2j-1}\over \sqrt j}\Big\}.   \eqno(3.4)  $$
  Then $k$ is well defined, $k\ge 2$, and 
$$ {N\o\over
 C_0     }    \le  {4^{2k-1}\over \sqrt k} \qq{\rm and} \qq {N\o\over
 C_0     }   >  {4^{2k-3}\over \sqrt{ k-1}}.$$
 Further 
$$k\le 3\log \big({N\o\over
 C_0     }\big).  \eqno(3.5 ) $$
Indeed, put for a while $X={N\o/ C_0     }$ and observe that $k\le 2^k$ and $(7k/2)-6\ge k/2$ when $k\ge 2$. Then 
$$X> {4^{2k-3}\over \sqrt{ k-1}}>{4^{2k-3}\over \sqrt{ k  }}\ge 2^{4k-6-k/2}=2^{(7k/2)-6 } \ge 2^{k/2}.$$
Hence $k\le (2/\log 2)\log X< 3\log X$. 
 \medskip By the assumption made,  the argument
$
\sum_{1\le
\ell\le N \, , j\not=
\ell}  \nu_\ell \l_\ell 
 $   is non-vanishing.
 Thus if $    \sup\{  |\nu_\ell|:  \ell\not=j\}>0$ 
$$\int_{d}^{d+T}e( t\sum_{1\le \ell\le
N \atop j\not=
\ell}  \l_\ell 
\nu_\ell)dt={e( d
  \sum_{1\le \ell\le
N \atop j\not=
\ell}  \l_\ell 
\nu_\ell   )\Big(e(  T
 \sum_{1\le \ell\le
N \atop j\not=
\ell}  \l_\ell 
\nu_\ell )-1\Big)\over 2i\pi\sum_{1\le \ell\le
N \atop j\not=
\ell}  \l_\ell 
\nu_\ell}, $$
and so 
$$\Big|\int_{d}^{d+T}e( t
  \sum_{1\le \ell\le
N \atop j\not=
\ell}  \l_\ell 
\nu_\ell   )dt\Big|=\Big| {  \sin  \pi T
  \sum_{1\le \ell\le
N \atop j\not=
\ell}  \l_\ell 
\nu_\ell \over  \pi \sum_{1\le \ell\le
N \atop j\not=
\ell}  \l_\ell 
\nu_\ell}\Big|. $$
  
Therefore 
$$\int_{d}^{d+T}{\bf \Upsilon}_j(t, \underline{  \b})dt = T\P\{S_k=0\}^{N-1}+H_j,$$
with$$H_j= \sum_{0< |\nu_\ell|\le (m-1)k\atop \ell\not=j} \Big(\prod_{1\le \ell\le
N \atop j\not=
\ell} \P\{S_k=\nu_\ell\} \Big)\int_{d}^{d+T}e\big( \sum_{1\le \ell\le
N \atop j\not=
\ell} (t\l_\ell-
\b_\ell)
\nu_\ell\big) dt,$$
 and
$$|H_j|\le     \sum_{0<  {\sup_{ \ell\not=j}}|\nu_\ell|\le (m-1)k  }    \Big(  \prod_{  \ell \not=j
  } \P\{Y_\ell=\nu_\ell\}    \Big) \Big| {  \sin  \pi T
  \sum_{1\le \ell\le
N \atop j\not=
\ell}  \l_\ell 
\nu_\ell \over  \pi \sum_{1\le \ell\le
N \atop j\not=
\ell}  \l_\ell 
\nu_\ell}\Big|     .$$
 
  It follows that 
$$\eqalign{ \int_{d}^{d+T}{\bf \Upsilon}(t, \underline{\b})dt    &\le   \big({\o\over
2m     }\big)^{2k} \Big( { N T \P\{S_k=0\}^{N-1} }  + 
\sum_{j=1}^N | H_j|\Big). \cr}\eqno(3.6) $$
Similarly 
    
$$  \int_{d}^{d+T}{\bf \Upsilon}(t, \underline{\b})dt= T\P\{S_k=0\}^{N } + H,$$
and 
$$|H |\le     \sum_{0<  {\sup_{ \ell }}|\nu_\ell|\le (m-1)k  }    \Big(  \prod_{  \ell 
  } \P\{Y_\ell=\nu_\ell\}    \Big) \Big| {  \sin  \pi T
  \sum_{1\le \ell\le
N }  \l_\ell 
\nu_\ell \over  \pi \sum_{1\le \ell\le
N  }  \l_\ell 
\nu_\ell}\Big|     .$$  So that
$$ T\P\{S_k=0\}^{N }-|H|\le  \big({\o\over
2m     }\big)^{2k} \big( {N   \over
   \P\{S_k=0\}    }\big)T\P\{S_k=0\}^{N }   + 
 \big({\o\over
2m     }\big)^{2k} \sum_{j=1}^N | H_j|.\eqno(3.7) $$  
  By the choice made in (3.4)  of $m$ and  $k$, we have
$$  \big({\o\over
2 m     }\big)^{2k}  N   
      = {N\over
4 ^{2k}    }        
       =   {N\o\over
2m4 ^{2k-1}    }        
       \le { C_0\over 2m\sqrt k}\le  {1\over 2} \P\{S_k=0\} .  \eqno(3.8) $$ 
  We get from (3.7) and (3.8)
$$ T\P\{S_k=0\}^{N }  \le 2\bigg(|H |     + 
\big({\o\over
2m     }\big)^{2k}    \sum_{j=1}^N | H_j|\bigg)\le 2\bigg(|H  |    + 
    {\P\{S_k=0\}\over 2N}   \sum_{j=1}^N | H_j|\bigg).   $$  
Hence 
$$ T\P\{S_k=0\}^{N }  \le 3\max\big (|H  |,\max  _{j=1}^N | H_j|\big).\eqno(3.9)  $$  

  We shall now bound  $|H_j|$ and $|H|$. We begin with $|H|$ and put ${\cal Z }_N=\sum_{1\le \ell\le
N  }   \l_\ell 
Y_\ell $.   We have 
$$|H |\le     \E\,   \Big| {  \sin  \pi T
  {\cal Z }_N \over  \pi   {\cal Z }_N}\Big|\cdot \chi\{ {\cal Z }_N\not= 0\}   .\eqno(3.10) $$
Indeed $$|H |\le     \sum_{0<  {\sup_{ \ell }}|\nu_\ell|\le (m-1)k  }    \Big(  \prod_{  \ell 
  } \P\{Y_\ell=\nu_\ell\}    \Big) \Big| {  \sin  \pi T
  \sum_{1\le \ell\le
N }  \l_\ell 
\nu_\ell \over  \pi \sum_{1\le \ell\le
N  }  \l_\ell 
\nu_\ell}\Big| = \E\,   \Big| {  \sin  \pi T
  {\cal Z }_N \over  \pi   {\cal Z }_N}\Big|\cdot \chi\{ {\cal Z }_N\not= 0\}   .  $$
 And we have the trivial bound
$$\E\,   \Big| {  \sin  \pi T
  {\cal Z }_N \over  \pi   {\cal Z }_N}\Big|\cdot \chi\{ {\cal Z }_N\not= 0\}\le\E\,     { 1 \over  \pi  |  {\cal Z }_N |} \cdot
\chi\{ {\cal Z }_N\not= 0\}\le {1
  \over  \pi  \displaystyle{\min_{0<  {\sup_{ \ell }}|\nu_\ell|\le (m-1)k  }\big\{|\sum_{1\le \ell\le
N  }   \l_\ell 
\nu_\ell |\big\}}   }. $$
By (3.4), (3.5), $mk= 2\o k\le  6\o  \log ( N\o/ C_0) $. Notice by using assumption (1.6) that  
$$\eqalign{   \min\big\{\big |\sum_{1\le \ell\le
N }    \l_\ell 
\nu_\ell \big|:0<   \sup_{ \ell } |\nu_\ell|\le (m-1)k     \big\}&=   \min_{{ u_k \,{\rm integers} \atop |u_\ell|\le 
(m-1)k}\atop  |u_1\l_1+\ldots+ u_N\l_N |\not=0 } \big |\sum_{1\le \ell\le N }    \l_\ell 
\nu_\ell \big|\cr &\ge   \min_{{ u_k \,{\rm integers} \atop |u_\ell|\le 
6\o  \log ( N\o/ C_0)}\atop  |u_1\l_1+\ldots+ u_N\l_N |\not=0 } \big |\sum_{1\le \ell\le N }    \l_\ell 
\nu_\ell \big| =\Xi\cr}$$
 Thus 
$$|H|\le {1\over \pi\Xi } .\eqno (3.11) $$
Similarly, letting ${\cal Z }_{N,j}=\sum_{1\le \ell\le
N\atop\ell\not=j  }   \l_\ell 
Y_\ell $, we have 
$$\eqalign{|H _j|&\le     \sum_{0<  {\sup_{ \ell\not=j}}|\nu_\ell|\le (m-1)k  }   \Big(  \prod_{  \ell 
  } \P\{Y_\ell=\nu_\ell\}    \Big) \Big| {  \sin  \pi T
  \sum_{1\le \ell\le
N\atop\ell\not=j }  \l_\ell 
\nu_\ell \over  \pi \sum_{1\le \ell\le
N \atop\ell\not=j }  \l_\ell 
\nu_\ell}\Big| \cr &= \E\,   \Big| {  \sin  \pi T
  {\cal Z }_{N,j} \over  \pi  {\cal Z }_{N,j}}\Big|\cdot \chi\{ {\cal Z }_{N,j}\not= 0\}   . \cr}$$
And so,
$$|H _j|\le \E\,   \Big| {  \sin  \pi T
   {\cal Z }_{N,j} \over  \pi    {\cal Z }_{N,j}}\Big|\cdot \chi\{  {\cal Z }_{N,j}\not= 0\} \le {1
  \over  \pi  \displaystyle{ \min_{0<  \sup_{ \ell\not=j}  |\nu_\ell|\le (m-1)k  }\big\{|\sum_{1\le \ell\le
N  }   \l_\ell 
\nu_\ell |\big\}}   }\le {1\over \pi \Xi}. \eqno (3.12)$$
By inserting these estimates into (3.9), we get 
$$ T\P\{S_k=0\}^{N }  \le  {3\over \pi \Xi}.\eqno(3.13)  $$


By (2.3),   
  $  \P\{S_k=0\}  \ge  C_0/(m\sqrt k)=C_0/(2\o\sqrt k)$, and by   reporting this into (3.13) and using (3.4), (3.5), we arrive to
 $$T    \le {3 \over \pi \Xi} \big( { 2\o\sqrt k \over   
      C_0} \big)^{N}\le  {3 \over \pi \Xi} \bigg( { 2\o\sqrt{3\log \big({N\o\over C_0}\big)  } \over   
      C_0} \bigg)^{N}. \eqno(3.14)$$
 
Consequently, if    
$$T    >  {3 \over \pi \Xi} \bigg( { 2\sqrt 3\,\o  \over   
      C_0} \sqrt{ \log \ {N\o\over C_0}   }\bigg)^{N},  $$
  then to     any reals $d,  
\b_1,\ldots, 
\b_N  $     corresponds a real $t\in [d,d+T]$  such that  
$$ \sup_{j=1}^N\, \lt  t\l_j-\b_j \rt\le {1\over \o}.  \eqno(3.15) $$ 
 The proof is now complete.
\cqfd


\bigskip   Theorem 1 has interesting consequences for  Dirichlet polynomials and more general polynomials. We shall investigate them.  Let
$  
\a_1,\ldots, 
\a_L  $ be given reals and consider the Dirichlet polynomials 
 $ D_L(t)= \sum_{n=1}^L \a_n n^{ it}$. 
Let $\pi(x)=\#\{p\, {\rm prime}\le x\}$ be  the prime number function. Choose $N=\pi(L)$. Using the prime factor   decomposition,  
  $n=p_1^{ a_1}\ldots p_N^{ a_N}$, $ a_j(n)\ge
0$,
$ 1\le j\le  N
$, $ 1\le n\le L$,   we get
$$ D_L(t)= \sum_{n=1}^L \a_n e^{ it\sum_{j=1}^N  a_j(n)\log p_j} .\eqno(3.16)$$
  Let $  
\t_1,\ldots, 
\t_N\in [0,1[$. Let $\O(n)=\sum_{j=1}^{ N }a_j(n)$ denotes the prime divisor function. As
$$\eqalign{ &\big|e^{2i\pi t\sum_{j=1}^N  a_j(n)\log
 p_j}-e^{2i\pi \sum_{j=1}^N  a_j(n)\t_j} \big|   \cr& =\big|e^{2i\pi\sum_{j=1}^N a_j(n)[(t\log  p_j-\nu_j -\t_j)+\nu_j
+\t_j]}-e^{2i\pi \sum_{j=1}^N  a_j(n)\t_j}\big|   =\big| e^{2i\pi
\sum_{j=1}^{ N }a_j(n) (t\log  p_j-\nu_j -\t_j)  }-1 
\big|  \cr&  \le 2 \pi
\sum_{j=1}^{ N }a_j(n) \big| t\log  p_j-\nu_j -\t_j\big| , \cr}  $$
by taking   the infimum over all $\nu_j$, we get
$$ \big|e^{2i\pi t\sum_{j=1}^N  a_j(n)\log
 p_j}-e^{2i\pi \sum_{j=1}^N  a_j(n)\t_j} \big|\le  2 \pi
\sum_{j=1}^  N  a_j(n) \lt t\log  p_j  -\t_j\rt    \le   2 \pi\big(\sup_{j=1}^  N \lt t\log  p_j  -\t_j\rt\big)  \O(n).$$  Herefrom   
$$\Big|  \sum_{n=1}^L \a_n n^{2i\pi t}- \sum_{n=1}^L \a_ne^{2i\pi \sum_{j=1}^N  a_j(n)\t_j} \Big|\le 
 2 \pi\big(\sup_{j=1}^  N \lt t\log  p_j  -\t_j\rt\big)   \cdot  \sum_{n=1}^L |\a_n|\O(n) . $$ 
Let $\o$ be some positive integer. By the comments made       after Theorem 1 concerning T\'uran's result, if $T    > T(N,\o):=  e^{   2 \o  N\log
{(N\o/ C_0)} 
\log N  }   $, then for any real $d$, any reals $\t_1,\ldots,\t_N$, there exists $\tau\in [d,d+T]$ such that 
$$ \sup_{j=1}^N\, \lt  \tau \log p_j-\t_j \rt\le {1/ \o}  .  $$ 
Thus
$$\Big|  \sum_{n=1}^L \a_n n^{2i\pi \tau}- \sum_{n=1}^L \a_ne^{2i\pi \sum_{j=1}^N  a_j(n)\t_j} \Big|\le 
{2 \pi\over \o} \cdot  \sum_{n=1}^L |\a_n|\O(n) . $$ 
 Let $\T=\R/\Z$ be the circle. Put for   $(\t_1,\ldots,\t_N)\in \T^{ N }$, $Q(\t_1,\ldots,\t_N)=\sum_{n=1}^L \a_ne^{2i\pi \sum_{j=1}^N 
a_j(n)\t_j}$.
\smallskip
 Consequently, given any
$(\t_1,\ldots,\t_N)\in \T^{ N }$,
$Q(\t_1,\ldots,\t_N)$ is well approached by $D_L(2 \pi \tau)$, for some $\tau\in [d,d+T]$ with an error term precised by the above
estimate. Now by (3.16), 
$$D_L(2 \pi \tau)=\sum_{n=1}^L \a_n e^{2i\pi \tau \sum_{j=1}^N  a_j(n)\log p_j} =Q(\lt\tau\log p_1\rt ,\ldots,\lt\tau\log
p_N\rt) .$$
Thereby
$$ 0\le  \sup_{(\t_1,\ldots,\t_N)\in \T^{ N }}\big| \sum_{n=1}^L
\a_ne^{2i\pi
\sum_{j=1}^N  a_j(n)\t_j}\big| -\sup_{  d   \le     \tau\le  d+T }\big|D_L(2 \pi \tau)\big| \le {2 \pi\over \o}\cdot  \sum_{n=1}^L
|\a_n|\O(n)  . \eqno(3.17)
$$
   Letting $-d$ and $T$ tend to infinity, next $\o$ tend to infinity yields (Bohr's reduction argument)
$$ \sup_{   t\in \R}\big|D_L(t)\big|=\sup_{(\t_1,\ldots,\t_N)\in \T^{ N }}\big| \sum_{n=1}^L \a_ne^{2i\pi
\sum_{j=1}^N  a_j(n)\t_j}\big|. 
$$
  Thus (3.17) means that 
$$ 0\le  \sup_{   t\in \R} |D_L(t) | -\sup_{ 2 \pi d   \le    t\le 2 \pi( d+T) } |D_L(t) | \le {2 \pi\over
\o}\cdot 
\sum_{n=1}^L |\a_n|\O(n)  , \eqno(3.18)
$$
 
Therefore the supremum of the Dirichlet polynomials $D_L$ over large intervals (of length greater than $T(N,\o)$) is
comparable to the supremum over the real line. And the error made is controlled by the degree of accuracy existing for the  Kronecker 
theorem within this interval. Further   estimate  (3.18) is {\it uniform} over $d$.\smallskip It would be interesting  to know below
which size of the interval this property breaks down.    Notice     by the
Dirichlet Theorem, that for any reals $  
\p_1,\ldots, 
\p_N  $  we
may choose $t\le \o^N$
such that 
$$ \sup_{j=1}^N\, \lt  t\p_j  \rt\le {1\over \o}.  \eqno(3.19)$$ 
(corresponding to the particular case  $  
\b_1=\ldots= 
\b_N=  d=0$ in Theorem 1). Further this is nearly optimal, see Erd\"os and R\'enyi's article [ER] for a discussion and for some related results
and the references therein, notably Haj\'os paper.  Therefore  this size cannot be smaller than $\o^N$.
\medskip
More generally, let $A$ be some positive real and  let
$\l_1,
\l_2,
\ldots,
\l_N 
$ be reals satisfying  
 $$\line{$ \qq \qq  \qq
\left\{  
 \matrix{ u_1\l_1+\ldots+u_N\l_N =0 \cr   
 \displaystyle{\max_{1\le j\le N}}
|u_j|\le 2A ,     \    \hbox{\rm $u_j\in \Z$}  \cr}  \right.  \qq \Longrightarrow\qq    
   u_1=u_2=\ldots=u_N=0.       \hfill$}   \eqno(3.20) $$
 Let $ a_j: \Z\to \Z $,
$ 1\le j\le  N
$ be arbitrary mappings,   and put
$$ {\cal N}_A=\Big\{b(n)=\sum_{j=1}^N  a_j(n)\l_j : \max_{1\le j\le N} 
|a_j(n)|\le A, n\in \Z \Big\}, \qq B(n)=\sum_{j=1}^N  a_j(n).  $$
  Because of assumption (3.20), to any $b\in {\cal N}_A$ corresponds a unique $n$ such that $b=b(n)$. Given $N$ reals $\a_1,\ldots,
\a_N$, consider   the   polynomials  
 $ {\cal D }_A(t)= \sum_{n   \in {\cal N}_A}\a_n e^{ itb(n)}$.  
  Let $  
\t_1,\ldots, 
\t_N\in [0,1[$.  Similarly
 $$\eqalign{ &\big|e^{2i\pi t\sum_{j=1}^N  a_j(n)\l_j}-e^{2i\pi \sum_{j=1}^N  a_j(n)\t_j} \big|    \le   2 \pi
\sum_{j=1}^  N  a_j(n) \lt t\l_j  -\t_j\rt    \le   2 \pi\big(\sup_{j=1}^  N \lt t\l_j  -\t_j\rt\big)B(n) .\cr} $$
And   $$\Big| \sum_{n   \in {\cal N}_A} \a_n n^{2i\pi t}- \sum_{n   \in {\cal N}_A} \a_ne^{2i\pi \sum_{j=1}^N  a_j(n)\t_j} \Big|\le 
 2 \pi\big(\sup_{j=1}^  N \lt t\l_j  -\t_j\rt\big)   \cdot  \sum_{n   \in {\cal N}_A} |\a_n|B(n) . $$ 
      Let $\o$  be  a positive integer  such that   $A< \o  \log {N\o\over C_0}$.  By   Theorem 1, if
$$T    >  {3 \over \pi \Xi(N,\o)} \bigg( { 2\sqrt 3\,\o  \over   
      C_0} \sqrt{ \log \ {N\o\over C_0}   }\bigg)^{N},  \qq  \qq \Xi(N,\o) = \min_{{ u_k \,{\rm integers} \atop |u_k|\le
  \o  \log {(N\o/ C_0)}}\atop  |u_1\l_1+\ldots+ u_N\l_N |\not=0 } \big |\sum_{1\le k\le N }    \l_k 
 u_k \big|, $$
  then to     any reals $d,  
\b_1,\ldots, 
\b_N  $     corresponds a real $\tau\in [d,d+T]$  such that    
 $ \sup_{j=1}^N\, \lt  \tau\l_j-\b_j \rt\le {1/\o}$.   
 \smallskip 
Consequently, by similar considerations
$$ 0\le \sup_{(\t_1,\ldots,\t_N)\in \T^{ N }}\big| \sum_{n   \in
{\cal N}_A}
\a_ne^{2i\pi
\sum_{j=1}^N  a_j(n)\t_j}\big|-\sup_{  d \le    \tau\le  d+T }\big|{\cal D}_A(2\pi \tau)\big|  \le {2 \pi\over
\o}\cdot 
\sum_{n   \in {\cal N}_A} |\a_n| B(n)  . \eqno(3.21)
$$
 Since, by letting $-d$, $T$,   next $\o$ tend to infinity 
$$ \sup_{   t\in \R}\big|{\cal D}_A(t)\big|=\sup_{(\t_1,\ldots,\t_N)\in \T^{ N }}\big|  \sum_{n   \in {\cal N}_A} \a_ne^{2i\pi
\sum_{j=1}^N  a_j(n)\t_j}\big|, 
$$
the same comments    concerning the supremums of the   polynomials ${\cal D}_A$ over large intervals are in order.


  \medskip\medskip
\noi   {\gum 4. Concludings remarks.} 
\medskip\noi We conclude this work by making several remarks related to the proof above and some key expressions having appeared in it,
as well as to some related questions.
\medskip\noi  
\item 1.  The central point of the proof is inequality (3.9):
$$ T\P\{S_k=0\}^{N }  \le 3\max\big (|H  |,\max  _{j=1}^N | H_j|\big).\eqno   $$
To get it, we had to adjust parameters $m$ and $k$   so that the factor $ \big({\o/
2m     }\big)^{2k}  ( {N /   \P\{S_k=0\}    } )$ of $T\P\{ S_k=0\}^N$ in
(3.7), can be made less than $1/2$. This operation seems inherent to the proof, thereby making the choice of $m$ and $k$ made in (3.4) 
unavoidable. Next
$|H|$ and
$|H_J|$ are controlled in exactly the same manner. For $H$  for instance, in (3.10)  we obtained the interesting bound 
$$|H |\le     \E\,   \Big| {  \sin  \pi T
  {\cal Z }_N \over  \pi   {\cal Z }_N}\Big|\cdot \chi\{ {\cal Z }_N\not= 0\}   ,  $$
and next continued with the rather brutal estimate $$\E\,   \Big| {  \sin  \pi T
  {\cal Z }_N \over  \pi   {\cal Z }_N}\Big|\cdot \chi\{ {\cal Z }_N\not= 0\}\le\E\,     { 1 \over  \pi  |  {\cal Z }_N |} \cdot
\chi\{ {\cal Z }_N\not= 0\}\le {1
  \over  \pi  \displaystyle{\min_{0<  {\sup_{ \ell }}|\nu_\ell|\le (m-1)k  }\big |\sum_{1\le \ell\le
N  }   \l_\ell 
\nu_\ell \big| }   }, $$
leading to  (3.11). At this stage, the question naturally arises whether this bound is really the best possible, in other
words how to compute 
$$  \E\,   \big| {  \sin  \pi T
  {\cal Z }_N \over  \pi   {\cal Z }_N}\big|  \chi\{ {\cal Z }_N\not= 0\}.\eqno(4.1 ) $$
 We believe that this is an important question.
  When in place of $ {\cal Z }_N$, we have a random variable $U$ with density distribution $G$, it is possible   to evaluate
 $\E\,  
\big| { 
\sin 
\pi T
 U\over  \pi  U}\big|$,  
by using the formula (see for instance   [K] p.430)
for any real $0<r<2$$$|x|^r={1\over 2K(r)}\int_{-\infty}^\infty { 1-\cos xt  \over |t|^{r+1}}\,dt={1\over K(r)}\int_{-\infty}^\infty { 
\sin^2 ({xt\over 2}) 
\over |t|^{r+1}}\,dt \eqno(4.2 )  $$ 
where $x$ is real and 
$$K(r)=  {\Gamma(2-r)\over
r(1-r)}\sin[(1-r){\pi\over 2}].\eqno(4.3 ) $$
    Choose   $1<r<2$.    By  writing
that
$|t|=|t|^{
  ({r+1\over 2}) }\cdot|t|^{  1-({r+1\over 2}) }$, next using  the
Cauchy-Schwarz inequality, we get by the aforementionned   formula  
$$\eqalign{  \E\,  
\big| { 
\sin 
x
 U\over   U}\big|=\int_\R \big| {\sin   {2xt\over 2})\over  t }\big|   G(t)dt  &\le \Big[\int_\R{\sin^2
{2xt\over 2}\over |t|^{r+1}}dt\Big]^{1/2}\cdot\Big[\int_\R{G^2(t)\over
|t|^{2[1-({r+1\over 2})]}} dt \Big]^{1/2}\cr &= \Big({|2x|^r \over 2K(r)}\Big)^{1/2}\cdot \Big[ \int_\R |t|^{r-1}G^2(t) 
   d t \Big]^{1/2}
,\cr}\eqno(4.4 )  $$
since $2[1-({r+1\over 2})]=2({1-r\over 2})=1-r $. Let $V$ be a random variable with density distribution $A_U^{-1}\cdot G^2(t)$
where  $A_U=   \int_\R G^2(t) dt $.
Thereby 
 $  \E\,  
\big| { 
\sin 
x
 U\over   U}\big|\le  \big({|2x|^r \over 2K(r)}  A_U \cdot  \E\,|V|^{r-1}\big)^{1/2} $. 
Letting $x=\pi T$, we obtain 
$$  \E\,  
\big| { 
\sin 
\pi T
 U\over   U}\big|\le  C_r \Big[ T ^{r } A_U \cdot  \E\,|V|^{r-1}\Big]^{1/2}.$$
 \medskip\noi 
 \item 2.   
 The construction made in Section 2 leads to an interesting observation concerning the general study of small
deviations in probability theory. The problem of evaluating 
$$ \P\{ |{\cal Z }_N|<\e\} $$
which is clearly related to the one of estimating  $\E\,   \big| {  \sin  \pi T
  {\cal Z }_N \over  \pi   {\cal Z }_N}\big| \, \chi\{ {\cal Z }_N\not= 0\}$, is of an arithmetic nature. And  so it seems that in
general, one cannot expect to find estimates of the small deviations of sums of i.i.d. random variables (even discrete and bounded) by
means on purely probabilist arguments only. The  intriguing remainding question is then to know which kind of conditions   on
the sequence $\l_n$, $n\le N$, would permit  to get sharp estimates of the small deviations.
 \medskip\noi 
\item 3. In a very recent work,  we obtained 	an estimate of integral   (4.1). The proof is rather delicate and   will be
published elsewhere. Although the bounds we found are sharp, there are unfortunately not sharp enough to be incorporated in the proof (section 3), and to
provide significant new results. But we showed that the integral in (4.1) appears in a rather wide context and obtained other applications.

\bigskip\noi  {\it Final note.} While writing down the paper, Chen [C2] (December 2007) informed us that his theorem 1 in [C1] can  also
provide another estimate for $t$, but different than ours and concerning $ \sum_{n=1}^N \lt t\l_n-\a_n\rt$. More precisely
 let $   \l_1, \l_2, \ldots ,\l_N $  be
linearly independent over $\Q$. Given $\e>0$,  
$$M_0=\Big[\big({N\pi^2\over8\e }\big)^{1/2}\Big]  ,\qq\quad \Lambda =
\min_{{ u_j \,{\rm integers} \atop |u_j|\le M_0}\atop  |u_1\l_1+\ldots+ u_N\l_N |\not=0 } \big|u_1\l_1+\ldots+ u_N\l_N\big|.$$

Put $$  T_0(\e, (\l_j))={NM_0^N\over 2\pi\Lambda }.$$
Let $     \a_1, \a_2, \ldots ,
\a_N
$ be real numbers. Then in any interval $J$ of length greater than $T_0(\e, (\l_j))$,  there exists a $t$ such that 
  $ \sum_{n=1}^N \lt t\l_n-\a_n\rt \le \e$.  Although the two quantities  $ \sum_{n=1}^N \lt t\l_n-\a_n\rt  $ and $
\sup_{n=1}^N \lt t\l_n-\a_n\rt $ are not really comparable, it is however interesting to compare the bounds for $T$ obtained in
each case, call them $T_C$ and $T_W$ respectively. Besides, Chen's approach and our are radically different.

\smallskip    {\it i)} Suppose we want to bound  $\sup_{n=1}^N \lt t\l_n-\a_n\rt $. Let $\e=1/\o$. Compare first $\Xi$ and $\Lambda$.
If $(\o/N)^{1/2}\log (N\o)={\cal O}(1) $, then $\Xi\gg \Lambda$. Next $   \log (T_C\Lambda)\sim     N\log (N\o) $ and $  \log (T_W\Xi)\sim
 (N\log
\o +\log\log (N\o)) $. Thus $T_W\ll T_C$. Now if $\o$ is large, namely if $(\o/N)^{1/2}\log (N\o)\not={\cal O}(1) $, then  $\Xi\ll
\Lambda$, the two preceding estimates of $T_C$ and $T_W$ remain valid, but   we do not see how to   compare
them.
   
\smallskip    {\it ii)} Suppose now we want to bound  $\sum_{n=1}^N \lt t\l_n-\a_n\rt $. Let $\e=\m^{-1}$, $\m$ integer and $\o=N\m$.
 Then  $   \log (T_W\Xi)\sim
 N\log
 (N\m) \sim \log (T_C\Lambda)  $.  The same comments     on $\Xi$ and $\Lambda$ are   in order.
\bigskip  \noi {\it Acknowledgments:} I wish to thank Mikha\"\i l Lifshits;   his    careful  reading of the preliminary versions, much helped to  
improve the writing of the
present paper.
  I
also thank Yong-Gao Chen for helpful comments.

\bigskip\noi{\gum References}
\medskip
    \noi [B]    Bacon H. M.  [1934]   {\sl An extension of Kronecker's Theorem},  Ann. of Math. {\bf 35},
776-786.
   \vskip 1pt \noi [C1]   Chen Y.-G.  [2000]   {\sl The best quantitative Kronecker's Theorem},  J. London Math. Soc. {\bf 61} no2,
691-701.
 \vskip 1pt \noi [C2]    Chen Y.-G.  [2007]   {\sl Private communication}.
  
\vskip 1pt \noi [ER]   Erd\H os P.,  R\'enyi A.  [1957] {\sl A probabilistic approach to problems of diophantine approximation}, 
Illinois J. Math.  {\bf 1}, 303-315.
  \vskip 1pt
\noi [K]    Kawata T.  [1972]  {\sl Fourier analysis in probability theory}, Academic Press, New York and London. 
  \vskip 1pt\noi [M] {\ph Montgomery H.} [1993]  {\sl  Ten lectures on the interface between analytic number theory and harmonic analysis},
Conference Board of the Math. Sciences, Regional Conference Series in Math. {\bf 84}.
    \vskip 1pt \noi [P]     Petrov V.  [1975]    {\sl Sums of independent
random variables}, Ergebnisse der Math. und ihre
Grenzgebiete, Springer  {\bf 82}. 
 \vskip 1pt \noi [Ti] Tijdeman R. [1989] {\sl Diophantine equations and diophantine approximations}, Number theory and applications
(Banff, AB, 1988), 215--243, NATO Adv. Sci. Inst. Ser. C Math. Phys. Sci  {\bf 265}, Kluwer Acad. Publ., Dordrecht, 1989.
 \vskip 1pt \noi [Tu] Tur\'an P. [1960] {\sl A theorem on diophantine approximation with application to Riemann zeta-function},  Acta
Math. Sci. Szeged   {\bf 21}, 311--318.
 \vskip 1pt \noi [W]   Weber M.  [2004]  {\sl Discrepancy of randomly sampled sequences of reals},  Math.
Nachrichten {\bf 271},   105-110.
   \medskip 
 \noi {\ph Address} : U.F.R. de Math\'ematique (IRMA), Universit\'e Louis-Pasteur et C.N.R.S., 7 rue Ren\'e
Descartes,  F-67084 Strasbourg Cedex  
 
\noi  {\ph Email} :  \tt  weber@math.u-strasbg.fr \par

\bye